\documentclass[12pt]{amsart}
\usepackage{amssymb,amsmath,amsthm,bm}
\usepackage[colorinlistoftodos,prependcaption,textsize=tiny]{todonotes}
 \definecolor{darkblue}{RGB}{0,0,160}
\usepackage{fouriernc} 
\usepackage[colorlinks=true,allcolors=darkblue]{hyperref}
\DeclareSymbolFont{usualmathcal}{OMS}{cmsy}{m}{n}
\DeclareSymbolFontAlphabet{\mathcal}{usualmathcal}

\hypersetup{
    colorlinks=true,
    linkcolor = blue,
    citecolor=magenta,
    urlcolor=cyan,
}

\usepackage[T1]{fontenc}

\usepackage{color}

\usepackage{parskip}
\usepackage{setspace}

\usepackage{amsmath,amsthm,amscd,amssymb, mathrsfs}

\usepackage{latexsym}

\usepackage{thm-restate}
\usepackage[capitalise, nameinlink]{cleveref}
\crefname{equation}{}{}
\crefname{figure}{{\sc Figure}}{{\sc Figure}}
\crefname{subsection}{Subsection}{Subsections}

\numberwithin{equation}{section}

\theoremstyle{plain}
\newtheorem{theorem}{Theorem}[section]
\newtheorem{lemma}[theorem]{Lemma}

\theoremstyle{definition}

\theoremstyle{remark}
\newtheorem{remark}[theorem]{Remark}

\newtheorem{case[theorem]}{Case}

\title[\parbox{14cm}{\centering{Parallelograms and VC dimensions\hspace{1in}}} \quad]{Parallelograms and the VC-dimension of the distance sets}
\author{Thang Pham}

\address{University of Science, Vietnam National University, Hanoi.}
\email{thangpham.math@vnu.edu.vn}

\begin{document}
\maketitle
\begin{abstract}
In this paper, we study the distribution of parallelograms and rhombi in a given set in the plane over arbitrary finite fields $\mathbb{F}_q^2$. As an application, we improve a recent result due to Fitzpatrick, Iosevich, McDonald, and Wyman (2021) on the Vapnik–Chervonenkis dimension of the induced distance graph. Our proofs are based on the discrete Fourier analysis. 
\end{abstract}
\section{Introduction}
Let $\mathbb{F}_q$ be a finite field of order $q$. Given $E_1, E_2, E_3, E_4\subset \mathbb{F}_q^2$, the set of parallelograms with one side-length $t$ in $E_1\times E_2\times E_3\times E_4$ is denoted by $Par_t(E_1, E_2, E_3, E_4)$, namely, 
\[Par_t(E_1, E_2, E_3, E_4):=\#\{(x, y, z, w)\in E_1\times E_2\times E_3\times E_4\colon x-y=z-w, ~||x-y||=t\}.\]
Here $||x-y||=(x_1-y_1)^2+(x_2-y_2)^2$.
Define
\[Par(E_1, E_2, E_3, E_4):=\sum_{t\in \mathbb{F}_q}Par_t(E_1, E_2, E_3, E_4),\]
as the set of all parallelograms. We also denote the set of rhombi with vertices in $E_1\times E_2\times E_3\times E_4$ of side-length $t$ by $Rhom_t(E_1, E_2, E_3, E_4)$.

The main purpose of this paper is to study the distribution of rhombi and parallelograms in $\mathbb{F}_q^2$ and an application on the problem of determining the Vapnik–Chervonenkis (VC) dimension of the induced distance graph in $\mathbb{F}_q^2$. For other similar configurations, say, rectangles and simplices, we refer the interested reader to \cite{alex, alex2, LM, FFA, shkredov, vinh} for more discussions.

Our first result is the following. 

\begin{theorem}\label{rhombus}
Given $t\ne 0$.
    Let $E_1, E_2, E_3, E_4\subset \mathbb{F}_q^2$. Assume that $|E_1||E_2|, |E_3||E_4|\gg q^3$, then
\[\left\vert Rhom_t(E_1, E_2, E_3, E_4)-\frac{Par_t(E_1, E_2, E_3, E_4)}{q}\right\vert\le q^{1/2}\left(\frac{|E_1||E_2|}{q}\right)^{1/2}\left(\frac{|E_3||E_4|}{q}\right)^{1/2}.\]
In addition, if $|E_1||E_2||E_3||E_4|\gg q^{7}$, then 
\[Rhom_t(E_1, E_2, E_3, E_4)\gg \frac{|E_1||E_2||E_3||E_4|}{q^4}.\]
\end{theorem}
This theorem offers a relation between $Rhom_t(E_1, E_2, E_3, E_4)$ and $Par_t(E_1, E_2, E_3, E_4)$. In this next theorem, we show that when $Par(E_1, E_2, E_3, E_4)$ is large enough, then 
\[Par_t(E_1, E_2, E_3, E_4)\sim \frac{Par(E_1, E_2, E_3, E_4)}{q}.\]
The formal statement reads as follows. 
\begin{theorem}\label{parallelgram}
Let $E_1, E_2, E_3, E_4\subset \mathbb{F}_q^2$ and $t\in \mathbb{F}_q^*$. We have 
\[\left\vert Par_t(E_1, E_2, E_3, E_4)-\frac{Par(E_1, E_2, E_3, E_4)}{q}\right\vert\le q^{1/2}(|E_1||E_2||E_3||E_4|)^{1/2}.\]
In addition, for any $t\ne t'$, one has 
\[\left\vert Par_t(E_1, E_2, E_3, E_4)-Par_{t'}(E_1, E_2, E_3, E_4) \right\vert \le q^{1/2}(|E_1||E_2||E_3||E_4|)^{1/2}.\]
\end{theorem}
We note that given sets $E_1, E_2, E_3, E_4$ in $\mathbb{F}_q^2$, the number parallelograms with vertices in $E_1\times E_2\times E_3\times E_4$ can be very large. A simple example is that 
\begin{equation}\label{ex-ad}E_1=E_2=E_3=E_4=\mathbb{F}_q\cdot \begin{pmatrix}1\\0\end{pmatrix}+X\cdot \begin{pmatrix}0\\1\end{pmatrix},\end{equation}
where $X$ is an arithmetic progression in $\mathbb{F}_q$. By a direct computation, the number of parallelograms with vertices in $E_1\times E_2\times E_3\times E_4$ is about $|X|^3q^3$. So as long as $|X|\gg q^{1/2}$, we have $Par_t(E_1, E_2, E_3, E_4)\sim |X|^3q^2$ and $Rhom_t(E_1, E_2, E_3, E_4)\sim |X|^3q$.

We now discuss an application of the two above theorems. Given $E\subset \mathbb{F}_q^2$ and $t\in \mathbb{F}_q\setminus \{0\}$. Define 
\[\mathcal{F}(E):=\{N(x)\colon x\in E\},\]
where $N(x)=\{y\in E\colon ||x-y||=t\}$.

We say that the system $\mathcal{F}(E)$ shatters a finite set $X\subset E$ if $\mathcal{F}(E)\cap X$ spans all subsets of $X$. We say the VC-dimension of $E$ is $d$ if there is a set $X\subset E$ of size $d$ such that it is shattered by $\mathcal{F}(E)$ and  no subset of size $d+1$ is shattered by $\mathcal{F}(E)$. 

The VC-dimension of the system $\mathcal{F}(E)$ has been studied recently by Fitzpatrick, Iosevich, McDonald, and Wyman \cite{ff-a}. More precisely, they proved the following theorem.

\begin{theorem}\label{f-a}
    Given $E\subset \mathbb{F}_q^2$ with $|E|\gg q^{\frac{15}{8}}$, then the VC-dimension of $\mathcal{F}(E)$ is equal to three.
\end{theorem}
In the next theorem, we present a threshold-improvement of this result, namely, we decrease the exponent $\frac{15}{8}$ to $\frac{13}{7}$. 

\begin{theorem}\label{mainthm}
Let $E\subset \mathbb{F}_q^2$. Assume that $|E|\gg q^{\frac{13}{7}}$, then the VC-dimension of $\mathcal{F}$ is equal to three. 
\end{theorem}
We do not know whether the exponent $13/7$ is sharp or not, and it is also difficult to suggest the next improvement. However, we do know that when the set $E$ has some nice structures, for instance, as in Example \ref{ex-ad}, a better threshold can be obtained. Some remarks about proofs of Theorems \ref{f-a} and \ref{mainthm} will be mentioned in Section \ref{last}.
\section{Preliminaries}
Let $f\colon \mathbb{F}_q^2\to \mathbb{C}$ be a complex-valued function. The Fourier transform of $f$ is defined by 
\[\widehat{f}(m)=\frac{1}{q^2}\sum_{x\in \mathbb{F}_q^2}f(x)\chi(-x\cdot m).\]
Here $\chi$ is a non-trivial additive character of $\mathbb{F}_q$. With this definition, the Fourier inversion formula and the Plancherel theorem are as follows:
$$ f(x)=\sum_{m\in \mathbb F_q^2} \chi(m\cdot x) \widehat{f}(m),~ \sum_{m\in \mathbb F_q^2} |\widehat{f}(m)|^2 =q^{-2}\sum_{x\in \mathbb F_q^2} |f(x)|^2.$$

Given a set $A\subset \mathbb{F}_q^2$, by abuse of notation, we denote its characteristic function by $A(x)$, namely, $A(x)=1$ if $x\in A$ and $0$ otherwise.

For any $j\ne 0$, let $S_j$ be the circle centered at the origin of radius $j$ defined as follows: 
$$S_j:=\{x\in \mathbb{F}_q^2\colon x_1^2+x_2^2=j\}.$$
The next lemma provides the precise form of the Fourier decay of $S_j$ for any $j\in \mathbb{F}_q$. 
A proof can be found in \cite{TAMS} or \cite{Kohsun}.
\begin{lemma}\label{large-fourier}
For any $j\in \mathbb{F}_q$, we have 
\[ \widehat{S_j}(m) = q^{-1} \delta_0(m) + \frac{1}{q^3} G_1^2(\eta, \chi) \sum_{r \in {\mathbb F}_q^*} 
\chi\Big(jr+ \frac{\|m\|}{4r}\Big),\]
where $\eta$ is the quadratic character, $\delta_0(m)=1$ if $m=(0,\ldots, 0)$ and $\delta_0(m)=0$ otherwise, and $G_1(\eta, \chi)$ is the Gauss sum.

In particular,  $|\widehat{S_j}(m)|\le q^{-3/2}$ if $m\ne (0, 0)$.
\end{lemma}
The following lemma is on the number of ``unit" distances in a pair of given sets in $\mathbb{F}_q^2$. A proof can be found in \cite{IR, Sh}.
\begin{lemma}\label{unit-distance}
Given $t\ne 0$.
    Let $U$ and $V$ be two sets in $\mathbb{F}_q^2$. We denote the number of pairs $(u, v)\in U\times V$ with $||u-v||=t$ by $N_t(U, V)$. Then we have 
\[\left\vert N_t(U, V)-\frac{|U||V|}{q}\right\vert\le 2q^{1/2}\sqrt{|U||V|}.\]
\end{lemma}

\section{Proof of Theorem \ref{rhombus}}
For $x$ with $||x||=t$, we denote $E_i\cap (E_j-x)$ by $E_{ij}^x$.
The strategy is as follows: if $u\in E_1\cap (E_2-x)$ and $v\in E_3\cap (E_4-x)$, then $u\in E_1, u+x\in E_2, v\in E_3, v+x\in E_4$. If $||u-v||=t$, then we have a rhombus of side-length $t$ with vertices $(u, u+x, v+x, v)$. 

Therefore, applying Lemma \ref{unit-distance}, one has
\[\left\vert Rhom_t(E_1, E_2, E_3, E_4)-\sum_{||x||=t}\frac{|E_{12}^x||E_{34}^x|}{q}\right\vert \le q^{1/2}\sum_{||x||=t}\sqrt{|E_{12}^x||E_{34}^x|}.\]

Assume that $|E_1||E_2|\gg q^3$ and $|E_3||E_4|\gg q^3$, then the second term is bounded from above by 
\[q^{1/2}\left(\frac{|E_1||E_2|}{q}\right)^{1/2}\left(\frac{|E_3||E_4|}{q}\right)^{1/2}.\]
The first term is exactly equal to $Par_t(E_1, E_2, E_3, E_4)/q$. Thus, the first estimate is proved. In the next step, we prove the lower bound for $Rhom_t(E_1, E_2, E_3, E_4)$. In particular, we want to show that 
\[Par_t(E_1, E_2, E_3, E_4)\gg \frac{|E_1||E_2||E_3||E_4|}{q^3}.\]
Since $|E_1||E_2|\gg q^3$ and $|E_3||E_4|\gg q^3$, the number of quadruples $(x, y, z, w)\in E_1\times E_2\times E_3\times E_4$ such that $||x-y||=||z-w||=t$ is $\gg \frac{|E_1||E_2||E_3||E_4|}{q^2}$. Thus, there is a direction with at least $|E_1||E_2||E_3||E_4|/q^3$ parallelograms with one side-length $t$. In other words, 
\[Par_t(E_1, E_2, E_3, E_4)\gg \frac{|E_1||E_2||E_3||E_4|}{q^3}.\]
This gives 
\[Rhom_t(E_1, E_2, E_3, E_4)\gg \frac{|E_1||E_2||E_3||E_4|}{q^4}-\frac{(|E_1||E_2||E_3||E_4|)^{1/2}}{q^{1/2}}\gg \frac{|E_1||E_2||E_3||E_4|}{q^4},\]
whenever $|E_1||E_2||E_3||E_4|\gg q^7$.
This completes the proof of the theorem. 
\begin{remark}\label{rm31}
In the above proof, when counting the number of quadruples $(u, u+x, v, v+x)\in E_1\times E_2\times E_3\times E_4$, it might happen that $v=u+x$. We refer to this tuple as a degenerate rhombus. The number of such rhombi is at most $\sum_{||x||=t}|E_{12}^x\ll |E_1||E_2|/q$. In practise, this scenario can be avoided completely if the sets $E_i$ are disjoint.  
\end{remark}
\section{Proof of Theorem \ref{parallelgram}}
Given $E_1, E_2, E_3, E_4\subset \mathbb{F}_q^2$. For any $t\ne 0$, we first show that $Par_t(E_1, E_2, E_3, E_4)$ is equal to
    \[q^6\sum_{m=-m'}\widehat{E_1}(-m)\widehat{E_2}(m)\widehat{E_3}(m)\widehat{E_4}(-m)\widehat{S_t}(0)+q^6\sum_{m\ne -m'}\widehat{E_1}(-m)\widehat{E_2}(-m')\widehat{E_3}(m)\widehat{E_4}(m')\widehat{S_t}(-m-m').\]
Indeed, 
\begin{align*}
&Par_t(E_1, E_2, E_3, E_4)=\sum_{v\in \mathbb{F}_q^2}\sum_{x, y\in \mathbb{F}_q^2}E_1(x)E_2(x+v)E_3(y)E_4(y+v)S_t(v)\\&=\sum_{x, y, v}\sum_{m, m'}E_1(x)E_2(y)\widehat{E_3}(m)\widehat{E_4}(m')\chi(m(x+v))\chi(m'(y+v))S_t(v)\\
&=q^6\sum_{m, m'}\widehat{E_1}(-m)\widehat{E_2}(-m')\widehat{E_3}(m)\widehat{E_4}(m')\widehat{S_t}(-m-m')\\
&=q^6\sum_{m=-m'}\widehat{E_1}(-m)\widehat{E_2}(m)\widehat{E_3}(m)\widehat{E_4}(-m)\widehat{S_t}(0)+q^6\sum_{m\ne -m'}\widehat{E_1}(-m)\widehat{E_2}(-m')\widehat{E_3}(m)\widehat{E_4}(m')\widehat{S_t}(-m-m')\\
&=I+II.
\end{align*}
To simplify the notations, we write $Par_t$ instead of $Par_t(E_1, E_2, E_3, E_4)$. 

For all non-zero $t$ and $t'$ with $t\ne t'$, we observe that $I(Par_t)=I(Par_{t'})$. On the other hand, using Lemma \ref{large-fourier}, we have
\begin{align*}
    &|II(Par_t)|=|q^6\sum_{m\ne -m'}\widehat{E_1}(-m)\widehat{E_2}(-m')\widehat{E_3}(m)\widehat{E_4}(m')\widehat{S_t}(-m-m')|\\
    &\le q^\frac{9}{2}\sum_{m\ne -m'}|\widehat{E_1}(-m)||\widehat{E_2}(-m')||\widehat{E_3}(m)||\widehat{E_4}(m')|\\
    &\le q^{9/2}\left(\sum_{m}|\widehat{E_1}(m)|^2\right)^{1/2}\left(\sum_{m}|\widehat{E_2}(m)|^2\right)^{1/2}\left(\sum_{m}|\widehat{E_3}(m)|^2\right)^{1/2}\left(\sum_{m}|\widehat{E_4}(m)|^2\right)^{1/2}\\
    &\le q^{1/2}(|E_1||E_2||E_3||E_4|)^{1/2}.
\end{align*}
This implies that \[|Par_t-Par_{t'}|\le q^{1/2}(|E_1||E_2||E_3||E_4|)^{1/2}.\]
As proved above,
 \[Par_t=q^6\left(\sum_{m=-m'}\widehat{E_1}(-m)\widehat{E_2}(m)\widehat{E_3}(m)\widehat{E_4}(-m)\widehat{S_t}(0)+\sum_{m\ne -m'}\widehat{E_1}(-m)\widehat{E_2}(-m')\widehat{E_3}(m)\widehat{E_4}(m')\widehat{S_t}(-m-m')\right).\]
The first term can be computed explicitly, namely,
\begin{align*}
I&=\frac{1}{q^3}\sum_{m}\sum_{x, y, z, w}E_1(x)E_2(y)E_3(z)E_4(w)\chi(m\cdot (x+y-z-w))=\frac{}{}
\frac{1}{q}\sum_{x, y, z, w}1_{x+y-z-w}=\frac{Par(E_1, E_2, E_3, E_4)}{q}.
\end{align*}
Since the absolute value of the second term is at most $q^{1/2}(|E_1||E_2||E_3||E_4|)^{1/2}$, we conclude that
\[\left\vert Par_t(E_1, E_2, E_3, E_4)-\frac{Par(E_1, E_2, E_3, E_4)}{q}\right\vert\le q^{1/2}(|E_1||E_2||E_3||E_4|)^{1/2}.\]
This completes the proof.
\section{Proof of Theorem \ref{mainthm}}\label{last}
To prove Theorem \ref{mainthm}, we need to find three distinct vertices $x^1, x^2, x^3\in E$ and vertices $y^{12}, y^{13}, y^{123}, y^1, y^2, y^3, y^0$ in $E$ such that the following hold:
 $$||y^{123}-x^1||=||y^{123}-x^2||=||y^{123}-x^3||=t.$$ 
 $$||y^{12}-x^1||=||y^{12}-x^2||=t, ~||y^{12}-x^3||\ne t.$$
 $$||y^{13}-x^1||=||y^{13}-x^3||=t, ~||y^{13}-x^2||\ne t.$$
 $$||y^{23}-x^2||=||y^{23}-x^3||=t, ~||y^{23}-x^1||\ne t.$$
   $$||y^{1}-x^1||=t, ~||y^{1}-x^2||\ne t, ~||y^{1}-x^3||\ne t.$$
  $$||y^{2}-x^1||\ne t, ~||y^{2}-x^2||= t, ~||y^{2}-x^3||\ne t.$$
   $$||y^{3}-x^1||\ne t, ~||y^{3}-x^2||\ne t, ~||y^{3}-x^3||=t.$$
 $$||y^{0}-x^1||\ne t, ~||y^{0}-x^2||\ne t, ~||y^{0}-x^3||\ne t.$$
This can be described clearly in the following picture.
\begin{center}
    \begin{figure}[h!]
\includegraphics[width=0.6\textwidth]{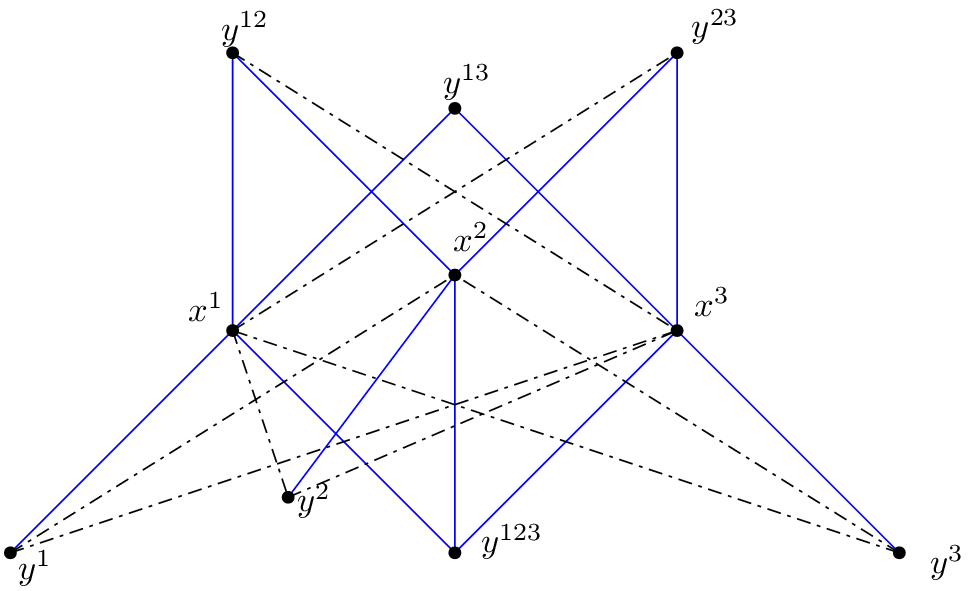}
\label{vcdim3figure}
\end{figure}
\end{center}
In this picture, if two vertices are of length $t$ then they are connected by a solid segment, otherwise by a dotted segment. 

Let $E'\subset E$ be the set of vertices $u$ such that the number of $v\in E$ and $||u-v||=t$ is at least $\frac{|E|}{2q}$. By Lemma \ref{unit-distance}, we have 
\[\frac{|E\setminus E'||E|}{q}-2q^{1/2}\sqrt{|E\setminus E'||E|}\le N(E\setminus E', E)<\frac{|E\setminus E'||E|}{2q}.\]
This implies 
\[|E\setminus E'|\le \frac{4q^3}{|E|^2}=o(|E|),\]
when $|E|\gg q^{13/7}$. A similar argument also shows that the set of vertices with at least $2|E|/q$ neighbors is $o(|E|)$. Thus, by abuse of notations, we might assume that each vertex in $E$ has at least $|E|/2q$ 
 and at most $2|E|/q$ neighbors in $E$. 

We now use Lemma \ref{unit-distance} again to find $u\in \mathbb{F}_q^2$ with $||u||=t$ such that $|E\cap (E-u)|\gg |E|^2/q^2$. Indeed, we have the number of pairs $(x, y)\in E^2$ such that $||x-y||=t$ is $(1-o(1))|E|^2/q$ and there are at most $q+1$ points on any circle of radius $t$, this implies the existence of such $u$. In the rest of the proof, we denote $E_u=E\cap (E-u)$, so $|E_u|\gg |E|^2/q^2$. Notice that $E_u\subset E$. 

The next step is to find a rhombus with vertices $(x^1, y^{123}, x^3, y^{13})\in E_u\times E_u\times E_u\times E$ such that $x^1-y^{123}\not\in\{ u, -u\}$, $y^{123}-y^{13}\not\in\{u, -u\}$, and $y^{123}-x^3\not\in \{u,-u\}$.

Let $A, B, C\subset E_u$ and $D\subset E$ be disjoint sets with $|A|=|B|=|C|\sim |E_u|$ and $|D|\sim |E|$. We want to find such a rhombus with vertices in $A\times B\times C\times D$. This helps us to avoid the degeneracy that might happen as noticed in Remark \ref{rm31}. 

Applying Theorem \ref{rhombus}, the number of rhombi with one side-length $t$ is at least 
\[\frac{|A||B||C||D|}{q^4},\]
under $|E|\gg q^{13/7}$.
Using the fact that the two circles of radius $t$ intersect in at most $2$ points, the number of rhombi $(x^1, y^{123}, x^3, y^{13})$ with $y^{123}-y^{13}\in \{u,-u\}$ is at most $2|E_u|$. The number of rhombi $(x^1, y^{123}, x^3, y^{13})$ with $y^{123}-x^1\in \{u,-u\}$ is at most $2N_t(E_u, E_u)$. With the same argument, the number of rhombi $(x^1, y^{123}, x^3, y^{13})$ with $y^{123}-x^3\in \{u,-u\}$ is at most $2N_t(E_u, E_u)$. Since $|E|\gg q^{13/7}$, $|E_u|\gg q^{12/7}$. Thus $N_t(E_u, E_u)=(1+o(1))|E_u|^2/q$. In total, the number of these rhombi is much smaller than $q^{-4}|A||B||C||D|$ whenever $|E|\gg q^{5/3}$. 

Hence, there are many rhombi $(x^1, y^{123}, x^3, y^{13})\in E_u\times E_u\times E_u\times E$ such that $x^1-y^{123}\not\in \{u, -u\}$, $y^{123}-y^{13}\not\in\{u, -u\}$, and $y^{123}-x^3\not\in \{u, -u\}$. We fix one of them. Set 
\[y^{12}=x^1+u, ~x^2=y^{123+u}, ~y^{23}=x^3+u.\]
Since $E_u=E\cap (E-u)$, we have $y^{12}, x^2, y^{23}\in E$. Since $y^{123}-x^1, y^{123}-x^3\not\in \{u, -u\}$, the three vertices $x^1, x^2, x^3$ are distinct. By our constraints, we can see that all vertices $x^1, x^2, x^3, y^{12}, y^{13}, y^{23}, y^{123}$ are distinct. 

The next step is to check $||y^{12}-x^3||\ne t$ and $||y^{23}-x^1||\ne t$. This is clear otherwise two distinct circles of radius $t$ intersect in more than two points. 

The final step is to choose $y^1, y^2, y^3$ with the desired property. This is also clear since we have assumed that each vertex of $E$ has at least $|E|/2q$ neighbors in $E$. Notice that $y^i\ne y^j$ for $i\ne j$, otherwise, we again have three points in the intersection of two circles of radius $t$.

Since each of vertices $x^1, x^2, x^3$ has at most $2|E|/q$ neighbors, the vertex $y_0$ can be chosen arbitrary in $E\setminus (N(x^1)\cup N(x^2)\cup N(x^3))$.

This completes the proof of the theorem.

\subsection{Discussions}
In this paragraph, we want to discuss briefly about the proof of Theorem \ref{mainthm} and that of Theorem \ref{f-a} in \cite{ff-a}. The main difference comes from the following step: in stead of counting the number of rhombi with vertices in $E_u\times E_u\times E_u\times E_u$ as in \cite{ff-a}, we replace one set $E_u$ by a bigger set $E$, this gives us more room to find desired configurations. Thus, the exponent $\frac{15}{8}$ can be decreased to $\frac{13}{7}$. In a recent work of Pham, Senger, Tait, and Thu-Huyen \cite{pseudo}, this VC-dimension notation has been studied in a more general setting for pseudo-random graphs. As mentioned in \cite{pseudo}, when applying their result to the distance graph, Theorem \ref{f-a} can not be improved. This can be explained in detail as follows. In \cite{pseudo}, to determine the VC-dimension, one has to bound the number of $C_4$ from below, but for pseudo-random graphs with the underlying vertex set $\mathbb{F}_q^2$, there is a possibility that it does not contain any non-trivial $C_4$. For example, the graph with the vertex set $\mathbb{F}_q^2$ and the two vertices $(x_1, x_2)$ and $(y_1, y_2)$ are connected by an edge if $x_1y_1+x_2y_2=1$. This graph is a pseudo-random graph. The adjacent relation can be viewed as the point $(x_1, x_2)$ belongs to the line defined by $y_1X+y_2Y=1$. Since two lines are either parallel or intersect in one point, we do not have any non-trivial $C_4$ in this graph, i.e. a $C_4$ with four distinct vertices. This explains why in $\mathbb{F}_q^2$, for the distance graph, the method in \cite{pseudo} is not applicable.  

If the reader is interested in the VC-dimension of other geometric configurations, we refer to some recent papers \cite{ABC, MSW, IMS} for more details. 
\section{Acknowledgements}
T. Pham would like to
thank to the Vietnam Institute for Advanced Study in Mathematics (VIASM) for the hospitality and for the excellent working conditions.

\bibliographystyle{amsplain}

\end{document}